\documentclass[12pt]{amsart}

\usepackage{amsmath}

\usepackage{amssymb}

\usepackage{enumitem}

\usepackage{graphicx,graphics,psfrag}
\usepackage{amscd}
\usepackage{amsthm}
\usepackage[mathscr]{eucal}
\usepackage{amssymb,latexsym,amsfonts}
\usepackage{fancyhdr}
\usepackage{float}
\usepackage[all]{xy}
\usepackage{graphicx}
\usepackage{framed}
\usepackage[colorlinks,linkcolor=blue, anchorcolor=blue,citecolor=blue]{hyperref}

\makeatletter
\@namedef{subjclassname@2020}{%
  \textup{2020} Mathematics Subject Classification}
\makeatother

\usepackage[T1]{fontenc}

\newcommand{\idealindex}{\mathrm{Index}}

\newtheorem{theorem}{Theorem}[section]
\newtheorem{corollary}[theorem]{Corollary}
\newtheorem{lemma}[theorem]{Lemma}
\newtheorem{proposition}[theorem]{Proposition}

\theoremstyle{definition}
\newtheorem{definition}[theorem]{Definition}
\newtheorem{remark}[theorem]{Remark}

\newtheorem{question}[theorem]{Question}

\numberwithin{equation}{section}

\begin{document}


\baselineskip=17pt


\title[The multiple points of maps from sphere to Euclidean space]{The multiple points of maps from sphere to Euclidean space}

\author[J. Wang]{Jun Wang}
\address{School of Mathematical Sciences\\
 Hebei Normal University\\
 Shijiazhuang\\
 Hebei, China}
\email{wjun@hebtu.edu.cn}

\author[X. Zhao]{Xuezhi Zhao}
\address{School of Mathematical Sciences\\
Capital Normal University\\
Beijing, China}
\email{zhaoxve@cnu.edu.cn}

\date{}

\begin{abstract}
In this paper, we obtain some sufficient conditions to guarantee the existence of multiple points of maps from $S^m$ to $\mathbb{R}^d$. Our main tool is the ideal-valued index of $G$-space defined by E. Fadell and S. Husseini.  	We obtain more detailed relative positional relationship of  multiple points. It is proved that for a continuous real value function  $f: S^m\rightarrow \mathbb{R}$  such that $f(-p)=-f(p)$, if $m+1$ is a power of $2$,  then there are  $m+1$  points $p_1, \ldots, p_{m+1}$ in $S^m$ such that $f(p_1)=\cdots=f(p_{m+1})$, where $p_1, \ldots, p_{m+1}$ are linearly dependent and any $m$ points of $p_1, \ldots, p_{m+1}$  are linearly independent. As a generalization of Hopf's theorem, we also  prove that for any continuous map $f: S^m\rightarrow \mathbb{R}^d$, if $m> d$, then there exists a pair of  mutually orthogonal points having the same image in addition to the antipodal points.
\end{abstract}



\maketitle

	\section{Introduction}\label{intro}
%
Many  interesting theorems have been proved about  how to guarantee the existence and relative positional relationship of multiple points of maps from sphere to Euclidean space.

One of famous problem was proposed by B. Knaster \cite{Knaster1947} in 1947: \par
\emph{Given a continuous map $f$ mapping the $(n-1)$-dimension Euclidean sphere $S^{n-1}\subset \mathbb{R}^{n}$ into $\mathbb{R}^m$, $m\leq n-1$, and $k=n-m+1$ points $p_1,\ldots ,p_k\in S^m$, does there exist a rotation $r\in SO(n)$ such that $f(r(p_1))=\cdots =f(r(p_k))$.} \par
 For the case $k=2$, which means $n-1=m$, H. Hopf's theorem \cite{hopf1944} answered this question.  More positive and counterexample are considered in \cite{Floyd1955, Marija2016}.

As  a special case of Knaster's problem, for the $m+1$-multiple points, Kakutani-Yamabe-Yujob\^{o}'s theorem \cite{Kakutani1942, yamabe}  proved that a continuous map on an $m$-sphere maps the ends of some $m+1$ mutually orthogonal unit vectors into a single value.
Dyson's theorem \cite{Dyson1951} tells us that a continuous real-valued function on a $2$-sphere maps the four end points of some pair of orthogonal diameters into a single value.  C.-T. Yang \cite{Yang1954,Yang1955} generalized the theorems of Borsuk-Ulam\cite{borsuk1933} , Kakutani-Yamabe-Yujob\^{o} and Dyson in 1950s, while D.G.Bourgin \cite{Bourgin1955} proved independently.   Yang-Bourgin theorem states that given any map $f: S^m\rightarrow \mathbb{R}^d$, the set $\{x \in S^m \mid  f(x)=f(-x)\}$ is of dimension $\geq m-d$.

 In this paper,  it is assumed $d\geq 1, n\geq 2$ throughout the entire article. We obtain  some sufficient conditions in Section \ref{proof} to guarantee the existence of $n$-multiple points under the continuous map $f$.

The  ideal-valued indices (see \cite{FadellHusseinni1988} or Section \ref{index} for definition) of  Stiefel manifold and generalized configuration spaces $\widetilde{W_{k,n}}(S^m)$ (see Section \ref{proof} for more details) have enormous effects in the proof of our results.

By the  ideal-valued indices of  Stiefel manifold, we prove Theorem \ref{main} which states that for a continuous map $f: S^m\rightarrow \mathbb{R}^d$  such that $f(-p)=-f(p)$, if
$$m\geq n, N>(n-1)d,$$
 there are  mutually orthogonal points $p_1, \ldots, p_n$ in $S^m$ such that $f(p_1)=\cdots=f(p_n)$ where $N=\min\{j\mid m-n+2\leq j\leq m+1, \binom{m+1}{j}\equiv 1 (\mod 2)\}$.

In particular, removing restrictions on mapping,
we obtain  Theorem \ref{thm:m>d} which states that  if  $m>d$,   for any continuous map $f: S^m\rightarrow \mathbb{R}^d$,  there  exists a pair of  orthogonal  points  having the same image  in addition to the pair of antipodal points obtained in Yang-Bourgin theorem \cite{Yang1954,Yang1955, Bourgin1955}.  This  can be also seen as a generalization of Hopf's theorem \cite{hopf1944} which  told us that if $m=d$, then  for any prescribed  positive number $\delta$, there exists a geodesic arc of length $\delta$ in $S^m$ with that the two end points have the same image under continuous map $f:S^m\rightarrow \mathbb{R}^d$. \par
Via the  ideal-valued indices generalized configuration spaces $\widetilde{W_{k,n}}(S^m)$,
Theorem \ref{thm:m:m+1} proves that for a continuous map $f: S^m\rightarrow \mathbb{R}$  such that $f(-p)=-f(p)$,  if  $$md< N',$$
where $N'=\min\{j\mid 2\leq j\leq m+1, \binom{m+1}{j}\equiv 1 (\mod 2)\}$, then there are  $m+1$  points $p_1, \ldots, p_{m+1}$ in $S^m$ such that $f(p_1)=\cdots=f(p_{m+1})$, where $p_1, \ldots, p_{m+1}$ are linearly dependent and any $m$ points of $p_1, \ldots, p_{m+1}$  are linearly independent.
It could be expectable that generalized configuration space will detect more relative positional relationship of multiple points, such as some solution to Knaster’s problem.

%
%



	Our paper is organized as follows.
In Section \ref{index}, it is  reviewed that the definition and some results about the ideal-valued index. In Section \ref{cgcs}, we shall calculate the ideal-valued indices of generalized configuration spaces  $\widetilde{W_{k,n}}(S^m)$.  Through ideal-valued indices of Stiefel manifold and $\widetilde{W_{k,n}}(S^m)$, we obtain some sufficient conditions in Section \ref{proof} to guarantee the existence and relative positional relationship of multiple points of maps from sphere to Euclidean space.

	\section{The ideal-valued  index}\label{index}
	
	In this section, it is  given a brief account of the definition of  ideal-valued index and its related properties.  More details and applications can be found in \cite{blagojevic2011,FadellHusseinni1988, Jaworowski1988}.\par
	
	\begin{definition}(E. Fadell and S. Husseini \cite{FadellHusseinni1988})
		Let $X$ be a paracompact Hausdorff space admitting an action of a compact Lie group $G$, and $R$ be a commutative ring. The ideal-valued index of  the space $X$, which is  denoted by $\idealindex^G(X; R)$, is defined to be the kernel
		$$
		\mbox{Ker} (c^*:H^*(BG; R)\rightarrow H^*(EG\times_G X; R))
		$$
		of the ring homomorphism $c^*$,	 where  $c: EG\times_G X\rightarrow BG$ is the map induced by the projective map $EG\times X\rightarrow EG$.
	\end{definition}
	\begin{remark}
		It is worth mentioning that the ideal-valued index of $X$ is related to the action of $G$ on $X$.
	\end{remark}
	If $f:X\rightarrow Y$ is a $G$-equivariant map, then we have the following commutative diagram
	$$ \xymatrix
	{
		H^{*}(BG; R)\ar[rr]^{\cong}\ar[d]_{{c_2}^*}&&  H^{*}(BG; R)\ar[d]^{{c_1}^*}\\
		H^{*}(EG\times_G Y; R)\ar[rr]^{f^*} &&  H^*(EG\times_G X; R).
	}
	$$
	Thus we obtain that
	\begin{proposition}\label{include}
		Let $X$ and $Y$ be two $G$-spaces, $f:X\rightarrow Y$ be a  $G$-equivariant map, then
		$$
		\idealindex^{G}(Y; R)\subset  \idealindex^{G}(X; R).
		$$
	\end{proposition}
	
	As a special case, let $G=\mathbb{Z}_2$, $R=\mathbb{Z}_2$, then $\idealindex^{\mathbb{Z}_2}(X; \mathbb{Z}_2)$, which is the ideal-valued index of $X$,  shall be an ideal of the cohomology ring
	$$
	H^*(BG;R)=H^*(\mathbb{R}P^{\infty};\mathbb{Z}_2)\cong\mathbb{Z}_2[\xi],\ \ \dim \xi=1.
	$$
	\begin{proposition}\label{Sn}
		Consider the $\mathbb{Z}_2$-action on the sphere $S^m$ which is induced by the antipodal map, then $\idealindex^{\mathbb{Z}_2}(S^m; \mathbb{Z}_2)$  is equal to  the ideal $(\xi^{m+1})$ in $\mathbb{Z}_2[\xi]$ generated by $\xi^{m+1}$.
	\end{proposition}
For the Stiefel manifold $V_{m+1,n}$ of all orthogonal $n$-frames in $\mathbb{R}^{m+1}$, Z.Z. Petrovi\'{c} \cite{zoran1997} consider the action of $\mathbb{Z}_2$, where the generator of $\mathbb{Z}_2$ acts by sending $(x_1,\cdots,x_n)$ to $(-x_1,\cdots,-x_n)$, and get the following results.
\begin{theorem}\cite[Theorem 3]{zoran1997} \label{thm:index:Vmn}
 $\idealindex^{\mathbb{Z}_2}(V_{m+1,n})$ is equal to the ideal $ (\xi^N)$ in $\mathbb{Z}_2[\xi]$ generated by $\xi^{N}$ where
 $$
 N=\min\{j\mid m-n+2\leq j\leq m+1, \binom{m+1}{j}\equiv 1 (\mod 2)\}.
 $$
\end{theorem}
From the above theorem, we can see that $m-n+2\leq N\leq m+1$. The following result is about when $N$ takes its maximal value and minimal value via the proof in \cite{zoran1997}.
\begin{proposition}\cite[Proposition 3]{zoran1997} \label{pro:N:maximal}
 For all $n\leq m$, $N=m+1$ if and only if $\nu_2(m+1)\geq \lceil \log_2(n)\rceil  $, where  $\nu_2(m+1)$ is the exponent of the highest power of $2$ dividing $m+1$, i.e. $\nu_2(m+1)=\max\{u\mid m+1 \equiv 0 \mod 2^u \}$.  \par
  Specially, if $n=m$, then $N=m+1$ if and only if $m+1$ is a power of $2$.
\end{proposition}
\begin{proposition}\label{pro:N:minimal}
   For all $1\leq n\leq m$, $N=m-n+2$ if and only if $m=2^\ell-2$ for some $\ell$.
\end{proposition}
	
		\section{The ideal-valued indices of generalized configuration spaces}\label{cgcs}
	
		Recall that the generalized configuration space of sphere \cite{Wang2019} is defined by
		\begin{align*}
			W_{k,n}(S^m)=\{&(p_1,\ldots ,p_n)\mid  p_i\in S^m, 1\leq i\leq n; \mbox{ for any k-elements subset }\\
			& \{i_1,\ldots ,i_{k}\}\subset \{1,\ldots ,n\},p_{i_1},\ldots ,p_{i_{k}} \mbox{ are linearly independent}\}.
		\end{align*}
		It is easy to see that
$$
W_{n,n}(S^m)\subset \cdots \subset W_{k+1,n}(S^m)\subset W_{k,n}(S^m)\subset W_{k-1,n}(S^m)\cdots \subset W_{1,n}(S^m).
$$
  Let $\widetilde{W_{k,n}}(S^m)=W_{k,n}(S^m)-W_{k+1,n}(S^m)$. Specially, when $k=n$, $\widetilde{W_{n,n}}(S^m)=W_{n,n}(S^m)$.  In this section, we shall calculate the ideal-valued indices of   generalized configuration spaces $\widetilde{W_{k,n}}(S^m)$, which will be used into the proof of our main result. \par
		Unless otherwise specified, the actions of $\mathbb{Z}_2$ on  configuration spaces $\widetilde{W_{k,n}}(S^m)$ and Stiefel manifolds $V_{m,n}$ are defined by
		\begin{align*}
			\mathbb{Z}_2\times \widetilde{W_{k,n}}(S^m) &\rightarrow \widetilde{W_{k,n}}(S^m)\\
			(\tau,x_1,\ldots ,x_n)&\mapsto (-x_1,\ldots ,-x_n),
		\end{align*}
		and
		\begin{align*}
			\mathbb{Z}_2\times V_{m,n}&\rightarrow V_{m,n}\\
			(\tau,y_1,\ldots ,y_n)&\mapsto (-y_1,\ldots ,-y_n),
		\end{align*}
		where $\tau$ is the generator of $\mathbb{Z}_2$.
		\par
		It is difficult to calculate the index directly because the cohomology of generalized configuration space $\widetilde{W_{k,n}}(S^m)$ is not known yet. But by the  relation between generalized configuration spaces $\widetilde{W_{k,n}}(S^m)$ and Stiefel manifolds, we can obtain the  following result.
	

		\begin{theorem}\label{thm:index:m:m+1}
			$\idealindex^{\mathbb{Z}_2}(\widetilde{W_{m,m+1}}(S^m);\mathbb{Z}_2)$ is equal to the ideal $(\xi^{N'}) $  in $\mathbb{Z}_2[\xi]$ generated by $\xi^{N'}$, where $N'=\min\{j\mid 2\leq j\leq m+1, \binom{m+1}{j}\equiv 1 (\mod 2)\}.$
		\end{theorem}
		
		\begin{proof}
For generalized configuration spaces $\widetilde{W_{m,m+1}}(S^m)$ and $\widetilde{W_{m,m}}(S^m)$, define  a map $f: \widetilde{W_{m,m+1}}(S^m)\rightarrow \widetilde{W_{m,m}}(S^m)$ as  $f(x_1,\ldots ,x_{m+1})=(x_1,\ldots ,x_m)$, and a map $g: \widetilde{W_{m,m}}(S^m)\rightarrow \widetilde{W_{m,m+1}}(S^m)$ as
			$g(x_1,\ldots ,x_{m})= (x_1,\ldots ,x_m,\frac{x_1+\cdots+x_m}{\| x_1+\cdots+x_m\|})$.\par
			It is not hard to verify  that the maps $f$ and $g$ both are $\mathbb{Z}_2$-equivariant.
			Thus by Proposition \ref{include}, we obtain that
			\begin{align*}
				\idealindex^{\mathbb{Z}_2}(\widetilde{W_{m,m+1}}(S^m);\mathbb{Z}_2)\subset \idealindex^{\mathbb{Z}_2}(\widetilde{W_{m,m}}(S^m);\mathbb{Z}_2),\\			\idealindex^{\mathbb{Z}_2}(\widetilde{W_{m,m}}(S^m);\mathbb{Z}_2)\subset\idealindex^{\mathbb{Z}_2}(\widetilde{W_{m,m+1}}(S^m);\mathbb{Z}_2) ,
			\end{align*}
			then	 $\idealindex^{\mathbb{Z}_2}(\widetilde{W_{m,m+1}}(S^m);\mathbb{Z}_2)=
			\idealindex^{\mathbb{Z}_2}(\widetilde{W_{m,m}}(S^m);\mathbb{Z}_2)$.

			For the ideal-valued indices of the generalized configuration spaces $\widetilde{W_{m,m}}(S^m)$, we have the relation
			$$
			\idealindex^{\mathbb{Z}_2}(\widetilde{W_{m,m}}(S^m);\mathbb{Z}_2)
=\idealindex^{\mathbb{Z}_2}(V_{m+1,m};\mathbb{Z}_2)
		$$
			By  the ideal-valued index in Theorem \ref{thm:index:Vmn}, we obtain the theorem.
		\end{proof}
By Proposition \ref{pro:N:maximal}, Proposition \ref{pro:N:minimal} and the proof in \cite{zoran1997}, we get the following corollaries.
\begin{corollary}\label{cor:m:m+1:=m+1}
			$\idealindex^{\mathbb{Z}_2}(\widetilde{W_{m,m+1}}(S^m);\mathbb{Z}_2)$ is equal to the ideal $(\xi^{m+1}) $  in $\mathbb{Z}_2[\xi]$ generated by $\xi^{m+1}$, if and only if  $m+1$ is a power of 2.
		\end{corollary}
\begin{corollary}\label{cor:m:m+1:=minimal}
		$\idealindex^{\mathbb{Z}_2}(\widetilde{W_{m,m+1}}(S^m);\mathbb{Z}_2)$ is equal to
\begin{enumerate}
\item $(\xi^{2}) \in \mathbb{Z}_2[\xi]$ if $m=2^\ell-2$ for some $\ell$,
  \item $(\xi^{2}) \in \mathbb{Z}_2[\xi]$ if $m\equiv 1,2 \mod 4$,
  \item $(\xi^{4}) \in \mathbb{Z}_2[\xi]$ if $m\equiv 3,4 \mod 8$.
\end{enumerate}
		\end{corollary}

\section{Main results}\label{proof}
In this section, we will prove some sufficient conditions to guarantee the existence of multiple points of maps from sphere to Euclidean space. The main tool is the ideal-valued indices of Stiefel manifold $V_{m,n}$ and generalized configuration space $\widetilde{W_{k,n}}(S^m)$. Further more, the relative positional relationship of multiple points will be obtained either.
\begin{theorem}\label{thm:m>d}
		Let $f: S^m\rightarrow \mathbb{R}^d$ be a continuous map. If $m>d$, then there exists  a pair of mutually orthogonal points $p_1,p_2\in S^m$  such that $f(p_1)=f(p_2)$.
	\end{theorem}

\begin{proof}
Consider the action  $\mathbb{Z}_2 \times V_{m+1,2}\rightarrow V_{m+1,2}$ defined by
$$(\tau, p_1, p_2)\mapsto (p_2, p_1),$$
where $\tau$ is the non-trivial element of $\mathbb{Z}_2$. In \cite[Theorem 4.11]{Basu2024}, it is proved that the ideal-valued index of $V_{m+1,2}$ under the action of permuting the orthogonal vectors is equal to the ideal  $(\xi^{m}) $  in $\mathbb{Z}_2[\xi]$ generated by $\xi^{m}$. \par
			For the continuous map $f: S^m\rightarrow \mathbb{R}^d$, define  $\widetilde{f}: V_{m+1,2}\rightarrow \mathbb{R}^{d}$ as  $\widetilde{f}(p_1,p_2)=f(p_1)-f(p_2)$. If there does not exist  $(p_1,p_2)\in V_{m+1,2}$ such that $f(p_1)=f(p_2)$, then $\widetilde{f}(V_{m+1,2})\subset \mathbb{R}^{d}-\{0\}$ and we can construct a map
			$$
			g=\phi\circ \widetilde{f}: V_{m+1,2}\rightarrow S^{d-1}
			$$
			where $\phi: \mathbb{R}^{d}-\{0\}\rightarrow S^{d-1}$ is defined by $\phi(x)=\frac{x}{\mid x\mid}$.
			Considering the action of  $\mathbb{Z}_2$ on $S^{d-1}$ induced by  the antipodal map, it is not hard  to verify that  the map $g$ is  $\mathbb{Z}_2$-equivariant.  Thus by  Proposition \ref{include}, we obtain that
			\begin{align}\label{belong}
				\idealindex^{\mathbb{Z}_2}(S^{d-1};\mathbb{Z}_2)\subset \idealindex^{\mathbb{Z}_2}(V_{m+1,2};\mathbb{Z}_2).
			\end{align}
			Then it is  obtained that $(\xi^{d})\subset (\xi^{m})$,
			thus $ m\leq d$,  and it is in contradiction to the assumption of the theorem. Hence the proof is completed.
		\end{proof}
Theorem \ref{thm:m>d} can be seen as  a generalization of Borsuk-Ulam theorem or Hopf's theorem \cite{hopf1944} where Hopf proved the case $m=d$. \par

\begin{theorem}\label{main}
Let $f: S^m\rightarrow \mathbb{R}^d$ be a continuous map such that $f(-p)=-f(p)$. Then for the  triples $(m,d,n)$ with
$$m\geq n, N>(n-1)d,$$
 there are  mutually orthogonal points $p_1, \ldots, p_n$ in $S^m$ such that $f(p_1)=\cdots=f(p_n)$. Here,  $N=\min\{j\mid m-n+2\leq j\leq m+1, \binom{m+1}{j}\equiv 1 (\mod 2)\}.$
\end{theorem}
\begin{proof}
For the continuous map $f: S^m\rightarrow \mathbb{R}^d$,  define  $\widetilde{f}: V_{m+1,n}\rightarrow \mathbb{R}^{nd}$ as  $$\widetilde{f}(p_1,\ldots ,p_n)=(f(p_1),\ldots ,f(p_n)).$$
			Suppose that there does not exist  $(p_1,\ldots ,p_n)\in V_{m+1,n}$ such that $f(p_1)=\cdots=f(p_n)$, then
			we can construct a map
			$$
			g=\frac{\varphi \circ \widetilde{f}}{|\varphi \circ \widetilde{f}|}: V_{m+1,n}\rightarrow S^{(n-1)d-1}.
			$$
			where $\varphi(x_1,\ldots ,x_n)=(x_1-\frac{x_1+\cdots +x_n}{n},\cdots,x_n-\frac{x_1+\cdots +x_n}{n})$  is defined to be the vector of differences between all variables and their mean.\par
			It is not hard  to verify that  the map $g$ is  $\mathbb{Z}_2$-equivariant since  $f(-p)=-f(p)$.  Thus by  Proposition \ref{include}, we obtain that
 \begin{align}\label{belong1}
				\idealindex^{\mathbb{Z}_2}(S^{(n-1)d-1};\mathbb{Z}_2)\subset \idealindex^{\mathbb{Z}_2}(V_{m+1,n};\mathbb{Z}_2).
			\end{align}
			For the ideal-valued index of  the sphere, we have that $$\idealindex^{\mathbb{Z}_2}(S^{(n-1)d-1};\mathbb{Z}_2)=(\xi^{(n-1)d})\subset \mathbb{Z}_2[\xi]$$
by Proposition~\ref{Sn}.

 Theorem \ref{thm:index:Vmn} tells us that $\idealindex^{\mathbb{Z}_2}(V_{m+1,n};\mathbb{Z}_2)=(\xi^{N})\subset \mathbb{Z}_2[\xi],$ where $N=\min\{j\mid m-n+2\leq j\leq m+1, \binom{m+1}{j}\equiv 1 (\mod 2)\}.$
			Then by Formula (\ref{belong1}),  it is  obtained that $(\xi^{(n-1)d})\subset (\xi^{N})$,
			thus $ N\leq (n-1)d$ and it is in contradiction to the assumption of the theorem.
		\end{proof}
The above theorem  gives a sufficient condition to ensure the existence of  the set $\mathcal{A}_{n}(f; S^m, \mathbb{R}^d)$.  At the same time, it is also obtained that  the relative positional relationship of the pre-image points of diagonal set of $(\mathbb{R}^d)^n$  because $(p_1, \ldots, p_n)$ belong to $V_{m+1,n}$ means $p_1, \ldots, p_n$ are mutually orthogonal.\par
In particular, considering the value of $N$  when it takes the maximal value, we obtain the following results by Proposition \ref{pro:N:maximal}.
\begin{lemma}\label{lem:point:index=m+1}
Let $f: S^m\rightarrow \mathbb{R}^d$ be a continuous map such that $f(-p)=-f(p)$. Then for the  triples $(m,d,n)$ with
$$ m\geq\max\{n, (n-1)d\},\nu_2(m+1)\geq \lceil \log_2(n)\rceil,  $$
 there are  mutually orthogonal points $p_1, \ldots, p_n$ in $S^m$ such that $f(p_1)=\cdots=f(p_n)$.  Here,  $\nu_2(m+1)$ is the exponent of the highest power of $2$ dividing $m+1$, i.e. $\nu_2(m+1)=\max\{u\mid m+1 \equiv 0 \mod 2^u \}$.
\end{lemma}
\begin{lemma}\label{lem:point:SO:index=m+1}
Let $f: S^m\rightarrow \mathbb{R}$ be a continuous map such that $f(-p)=-f(p)$. Then if $m$ is a power of $2$,
 then there are  mutually orthogonal points $p_1, \ldots, p_{m+1}$ in $S^m$ such that $f(p_1)=\cdots=f(p_{m+1})$.
\end{lemma}
For any value of $m,n$,  the number $m-n+2$ is the minimal value of $N$ by Proposition \ref{pro:N:minimal}, then we obtain the following by Theorem \ref{main}.
\begin{lemma}\label{lem:point:index=m-n+2}
Let $f: S^m\rightarrow \mathbb{R}^d$ be a continuous map such that $f(-p)=-f(p)$. If
$$m\geq (n-1)(d+1),$$
 then there are  mutually orthogonal points $p_1, \ldots, p_n$ in $S^m$ such that $f(p_1)=\cdots=f(p_n)$.
\end{lemma}

	Given different values of $\{m, d, n\}$, we get some interesting results from above lemmas. In particular, let $n=2$, we get  the  following corollary by Lemma \ref{lem:point:index=m+1}.\par
	
	\begin{corollary}\label{cor:m:odd}
		Let $f: S^m\rightarrow \mathbb{R}^d$ be a continuous  map such that $f(-p)=-f(p)$.
		If $m$ is odd, $m\geq \max\{2,d\}$, then there is a pair of mutually orthogonal points $p_1,p_2\in S^m$ such that $f(p_1)=f(p_2)$.
	\end{corollary}
	
%
%

	Applying Lemma \ref{lem:point:index=m-n+2} with $n=2$,  we obtain the following corollary which is a similar result as in Corollary \ref{cor:m:odd}.

	\begin{corollary}\label{BU2}
		Let $f: S^m\rightarrow \mathbb{R}^d$ be a continuous map such that $f(-p)=-f(p)$.
		If  $m\geq d+1$,  then there exist
		 a pair of mutually orthogonal points $p_1,p_2\in S^m$ such that $f(p_1)=f(p_2)$.
	\end{corollary}
%
%

	By the ideal-valued index of generalized configuration spaces $\widetilde{W_{k,n}}(S^m)$, we can obtain more relative positional relationships of multiple points. Through the ideal-valued indices of  $W_{m, m+1}(S^m)$ (see Theorem \ref{thm:index:m:m+1}), we  obtain the following Theorem.
\begin{theorem}\label{thm:m:m+1}
		Let $f: S^m\rightarrow \mathbb{R}^d$ be a continuous map such that $f(-p)=-f(p)$. If  $$md< N',$$
then there are  $m+1$  points $p_1, \ldots, p_{m+1}$ in $S^m$ such that $f(p_1)=\cdots=f(p_{m+1})$, where $p_1, \ldots, p_{m+1}$ are linearly dependent and any $m$ points of $p_1, \ldots, p_{m+1}$  are linearly independent. Here,
$N'=\min\{j\mid 2\leq j\leq m+1, \binom{m+1}{j}\equiv 1 (\mod 2)\}$,
\end{theorem}
\begin{proof}
By the analogous arguments as in the proof of Theorem \ref{main}, if there is no   point $(p_1, \ldots, p_{m+1})\in \widetilde{W_{m,m+1}}(S^m)$  such that $f(p_1)=\cdots=f(p_{m+1})$,  we can construct a $\mathbb{Z}_2$-equivariant  map
			$$
			g=\frac{\varphi \circ \widetilde{f}}{|\varphi \circ \widetilde{f}|}: \widetilde{W_{m,m+1}}(S^m) \rightarrow S^{((m+1)-1)d-1}.
			$$
Then $\idealindex^{\mathbb{Z}_2}(S^{((m+1)-1)d-1};\mathbb{Z}_2)\subset \idealindex^{\mathbb{Z}_2}(\widetilde{W_{m,m+1}}(S^m);\mathbb{Z}_2).$

Due to Theorem \ref{thm:index:m:m+1}, and Proposition \ref{include},
$$
(\xi^{md}) =(\xi^{((m+1)-1)d}) \subset (\xi^{N'}) ,
$$
where $N'=\min\{j\mid 2\leq j\leq m+1, \binom{m+1}{j}\equiv 1 (\mod 2)\}$.
We obtain $md\geq N'$,  it is a contradiction.
\end{proof}
In particular, let $d=1$, we get the following.
\begin{corollary}\label{cor:m:m+1}
		Let $f: S^m\rightarrow \mathbb{R}$ be a continuous map such that $f(-p)=-f(p)$. If $m+1$ is a power of $2$,
then there are  $m+1$  points $p_1, \ldots, p_{m+1}$ in $S^m$ such that $f(p_1)=\cdots=f(p_{m+1})$, where $p_1, \ldots, p_{m+1}$ are linearly dependent and any $m$ points of $p_1, \ldots, p_{m+1}$  are linearly independent.
\end{corollary}
\begin{proof}
  If $m+1$ is a power of $2$,  by corollary \ref{cor:m:m+1:=m+1},  we get the ideal-valued index
  $$ \idealindex^{\mathbb{Z}_2}(\widetilde{W_{m,m+1}}(S^m);\mathbb{Z}_2)=(\xi^{m+1}).$$
Then we complete the proof via Theorem \ref{thm:m:m+1}.
\end{proof}
Corollary \ref{cor:m:m+1} can be seen as a generalization of Kukutani-Yamabe-Yujob\^{o}'s theorem \cite{Kakutani1942, yamabe}  which states that a continuous real-valued function on an $n$-sphere maps the terminals of some $n+1$ mutually orthogonal radii into a single value. The relative position relationship of multiple points in Corollary \ref{cor:m:m+1} is more precise.

When $N'$ takes the minimal value $2$,  we obtain
\begin{corollary}\label{cor:m:m+1:=2}
		Let $f: S^1\rightarrow \mathbb{R}^1$ be a continuous map such that $f(-p)=-f(p)$,
then there are  $2$  points $p_1, p_{2}$ in $S^1$ such that $f(p_1)=f(p_{2})$, where $p_1,  p_{2}$ are linearly dependent.
\end{corollary}
Corollary \ref{cor:m:m+1:=2} is actually a special case of Borsuk-Ulam theorem.
Similarly, by Corollary \ref{cor:m:m+1:=minimal}, if $m=3$, then we obtain the following corollary.
\begin{corollary}\label{cor:m:m+1:=2}
		Let $f: S^3\rightarrow \mathbb{R}^1$ be a continuous map such that $f(-p)=-f(p)$,
then there are  $4$  points $p_1, \ldots, p_{4}$ in $S^m$ such that $f(p_1)=\cdots=f(p_{4})$, where $p_1, \ldots, p_{4}$ are linearly dependent and any $3$ points of $p_1, \ldots, p_{4}$  are linearly independent.
\end{corollary}

By the definition of generalized configuration spaces, it is understood $W_ {k+1,n}(S^m)\subset W_{k,n}(S^m)$, then above discussion
suggests an  interesting question as follows.

\begin{question}
For any continuous map $f: S^m\rightarrow \mathbb{R}^d$, does there exist an element $(p_1,\ldots , p_{n})\in W_{k, n}(S^m)\backslash W_{k+1, n}(S^m)$ such that $f(p_1)=\cdots=f(p_{n})$?
\end{question}
Theorem \ref{thm:m:m+1} answers this question for special cases.
		

\section*{Acknowledgements}
The authors express their thanks to Professor Wac{\l}aw Marzantowicz for his helpful comments and suggestions.


\end{document}